\documentclass{amsart}
\usepackage{graphicx}
\usepackage{amsmath}
\usepackage{amsthm}
\usepackage[pdftex,colorlinks]{hyperref}
\usepackage[all]{xy}

\theoremstyle{definition} 
\theoremstyle{plain} \newtheorem{theorem}[subsection]{Theorem}
\theoremstyle{plain} 
\theoremstyle{plain} \newtheorem{proposition}[subsection]{Proposition}
\theoremstyle{plain} \newtheorem{lemma}[subsection]{Lemma}
\theoremstyle{plain} \newtheorem*{theorem*}{Theorem}

\theoremstyle{remark} 
\theoremstyle{remark} 

\mathchardef\csum="2023

\def\ksm{\mathsf{ ksm}}

\def\rok{\mathsf{ rok}}

\def\css{\mathsf{ css}}

\def\int{\mathsf{ int}}

\def\Z{\mathbf{ Z}}
\def\R{\mathbf{ R}}
\def\cone{\mathsf{ cone}}

\def\link{\mathsf{ link}}

\begin{document}
\title{The triangulation of manifolds}
\author{Frank Quinn}
\date{November 2013}
\thanks{This work was partially supported by the Max Planck Institute for Mathematics in Bonn.}
\begin{abstract}
A mostly expository account of old questions about the relationship between polyhedra and topological manifolds. Topics are old topological results, new gauge theory results (with speculations about next directions), and history of the questions. MR classifications  57Q15, 01A60, 57R58.
\end{abstract}
\maketitle

\tableofcontents

 \section{Introduction}
 This is a survey of the current state of triangulation questions posed by   
Kneser in 1924:
\begin{enumerate}\item Is a polyhedron with the local homology properties of Euclidean space, locally homeomorphic to Euclidean space? 
\item Is a space locally homeomorphic to Euclidean space, triangulable (homeomorphic to some polyhedron)?
\item If there are two such triangulations, must they be PL equivalent?
\end{enumerate}
 Topological work on the topic is  described in Sections  \ref{sect:nearMflds}-\ref{sect:cssinvariant}. This work was mature and essentially complete by 1980, but leaves open questions about \(H\)-cobordism classes of homology 3-spheres.  Gauge theory has had some success with these, with the most substantial progress for the trianguation questions made in a recent paper of Manolescu \cite{mano13}.  Manolescu's paper is discussed in Section \ref{sect:gaugetheory}. This area is not yet mature, and one objective is to suggest other perspectives.
 Section \ref{sect:history} recounts some of the history of Kneser's questions. Kneser posed them as an attempt to provide foundations for Poincar\'e's insights twenty years before. They were found to be a dead end without significant applications, but were  fruitful challenges to technology as it developed. The progress of the subject can be traced out in applications to these questions, but here we have a different  concern: why did Kneser point  his contemporaries into a dead end? Or was he trying to get them to face the fact that it was a dead end? The answers give a window into the transition from pre-modern to modern mathematics in the early twentieth century.

The remainder of the introduction gives modern context for the questions and describes the organization of the technical parts of the paper. 
 \subsection{Modern context}\label{sect:questions} 
 The relevant main-line topics are PL manifolds, topological manifolds, and ANR homology manifolds\footnote{ANR = `Absolute Neighborhood Retract'. For finite-dimensional spaces this is equivalent to `locally contractible', and is used to rule out local point-set pathology.}. Polyhedra that are homology manifolds, referred to here as ``PL homology manifolds'', are mixed-category objects, and Kneser's questions amount to asking how these are related to the main-line categories.
 \[\xymatrix{\text{PL manifold}\ar[r]^{\subset}\ar[dr]^{\subset}\ar[dr]^{\subset}&\text{Topological}\ar[r]^{\subset}&\text{ANR homology}\\
&\text{PL homology}\ar[ur]^{\subset}\ar@{<-->}[u]_{?}
}\]
The standard  categories differ radically in flavor and technique, but turn out to be almost the same. For the purposes here, topological and ANR homology manifolds are the same. PL and topological manifolds differ by the Kirby-Siebenmann invariant \(\ksm(M)\in H^4(M;Z/2)\). This is in a single cohomology group, with the smallest possible coefficients, so is about as small as an obstruction can be without actually being zero. This means the \emph{image} of PL homology manifolds in the main-line picture is highly constrained, and on the image level the answers to the questions can't be much different from `yes'.  Unfortunately it turns out that there are a great many PL homology manifolds in each image equivalence class.

\subsection{Topology}
PL homology manifolds have two types of singularities: dimension 0 (problematic links of vertices), and codimension 4 (problematic links of \((n-4)\)-simplices). 4-manifolds are special, in part because these two types coincide. Vertex singularities can be canonically resolved, so are topologically inessential. This is described in \S\ref{sect:nearMflds} and provides an answer to Kneser's first question. 

The real difficulties come from the codimension 4 singularities, and these involve homology spheres. 
We denote the group of homology H-cobordism classes of homology 3-spheres by \(\Theta\). The full official  name is \(\Theta^H_3\), but the decorations are omitted here because they don't change. Several descriptions of this group are given in \S\ref{ssect:theta3}. 

A PL homology manifold  \(K\) has an easily-defined Cohen-Saito-Sullivan cohomology class \(\css(K)\in H^4(K;\Theta)\) \cite{sullSing}, \cite{cohen}; see \S\ref{ssect:css}. The Rokhlin homomorphism \(\rok\colon \Theta\to Z/2\) induces a change-of-coefficients exact sequence
 \[\xymatrix{\ar[r]& H^4(K;\ker(\rok))\ar[r] &{H^4(K;\Theta)}\ar[r]^{\rok}& H^4(K;Z/2)\ar[r]^{\beta}& H^5(K;\ker(\rok))}\]
 with Bokstein connecting homomorphism \(\beta\). 
 The image of the Cohen-Saito-Sullivan class is the Kirby-Siebenmann class. The baby version of the main theorem is 
 \begin{theorem*}\cite{gs} If \(M\) is a topological manifold of dimension \(\neq4\) (and boundary of dimension \(\neq4\) if it is nonempty) then concordance classes of homeomorphisms to polyhedra correspond to lifts of the Kirby-Siebenmann class to \(H^4(M;\Theta)\).
 \end{theorem*}
 This is the baby version because serious applications (if there were any) would require the relative version, Theorem \ref{thm:main}.
A corollary is that a triangulation exists if and only if the Bokstein of the Kirby-Siebenmann invariant is trivial. Note that if the Kirby-Siebenmann class lifts to a class with integer coefficients then it lifts to any coefficient group, and it follows that the manifold is triangulable. Similarly, if the Kirby-Siebenmann class lifts to coefficients \(Z/k\) but no further, then triangulability of the manifold depends on whether or not there is an element  in \(\Theta\) of order \(k\) and nontrivial Rokhlin invariant. Finally, triangulations are classified up to concordance by \(H^4(K;\ker(\rok))\).

These results reduce the geometric questions to questions about the group \(\Theta\) and the Rokhlin homomorphism. This part of the picture was essentially complete by 1980, but \(\Theta\) is opaque to traditional topological methods.  It has grudgingly yielded some of its secrets to sophisticated gauge theory; an overview is given in \S\ref{sect:gaugetheory}. It is infinitely-generated and lots of these generators have infinite order. This means if there is a triangulation of \(M\) and \(H^4(M;Z)\neq0\) then there are are a great many different ones. Manolescu's recent advance is  that the Rokhlin homomorphism does not split. This implies that there are manifolds (eg.~the ones identified by Galewski-Stern \cite{gs5}) that cannot be triangulated. Manolescu's paper is described in \S \ref{sect:gaugetheory}. This theory is in a relatively early stage of development so the section gives speculations about future directions.

It  seems reasonable to speculate that homology spheres with nontrivial Rokhlin invariant must have infinite order. Indeed, it seems reasonable to expect \(\Theta\) to be torsion-free. Either would imply that that \(M\) has a triangulation if and only if the integral Bokstein \(\beta\colon H^4(M;Z/2)\to  H^5(M;Z)\) is trivial on the Kirby-Siebenmann class. Proof of existence of triangulations in such cases should be easier than cases that might involve torsion. 
 
\section{Homology manifolds are essentially manifolds}\label{sect:nearMflds} We begin with Kneser's first question because the answer is easy (now) and sets the stage for the others. 

\subsection{Singular vertices} Suppose \(L\) is a PL homology \(n\)-manifold with homology isomorphic to the homology of the \(n\)-sphere. Then the cone on \(L\) is a PL homology \((n+1)\)-manifold (with boundary). However if the dimension is greater than 1 (to exclude circles) and  \(L\) is \emph{not simply-connected} then the cone point is not a manifold point. The reason is that the relative homotopy group 
\[\pi_2(\cone , \cone - *)\simeq \pi_1(L)\neq \{1\}\]
is nontrivial, and this is impossible for a point in a manifold. 

The lowest dimension in which non-simply-connected homology spheres occur is 3, and the oldest and most famous 3-dimensional example was described by Poincar\'e, see Kirby-Scharlemann \cite{kscharlemann79}. There are examples in all higher dimensions but the 3-dimensional ones are the most problematic.  
These cone points turn out to be the only topological singularities:
\begin{theorem}\label{almosttop} A PL homology manifold is a topological manifold except at vertices with  non-simply-connected links of dimension greater than 2.\end{theorem}
This statement is for manifolds without boundary, but extends easily. Boundary point are singular if either the link in the boundary, or the link in the whole manifold, is non-simply-connected.
 
We give a quick proof using mature tools from the study of ANR homology manifolds. Most homology manifolds are not manifolds, and some of them are quite ghastly\footnote{There is a technical definition of `ghastly' in \cite{davermanwalsh} that lives up to the name.}. Nonetheless they are close to being manifolds. There is a single obstruction in \(H^0(X;Z)\) whose vanishing corresponds to the existence of a map \(M\to X\) with essentially-contractible point inverses, and \(M\) a topological manifold \cite{Q1}. These are called \emph{resolutions} by analogy with resolution of singularities in algebraic geometry. When a resolution exists it is unique, essentially up to homeomorphism. Roughly speaking this gives an equivalence of categories, and the global theories are the same. 

The obstruction is so robust that a heroic effort was required to show that exotic examples exist \cite{BFMW}. Existence of a manifold point implies the obstruction vanishes, so PL homology manifolds have resolutions.

Next, Edwards' CE approximation theorem asserts that if \(X\) is an ANR homology manifold of dimension at least 5, and \(r\colon M\to X\) is a resolution, then \(r\) can be approximated by a homemorphism if and only if \(X\) has the `disjoint 2-disk property', see \cite{daverman}. It is easy to see that PL homology manifolds of dimension at least 5 have the disjoint disk property everywhere except at \(\pi_1\)-bad vertices. This completes the proof except in dimension 4, where the only question is with cones on homotopy spheres. Perelman has shown these are actually standard, so the cone is a PL 4-ball and the cone point is a PL manifold point. The weaker assertion that they are topologically standard also follows from the next section.

This proof seems effortless because we are using big hammers on small nails. The job could be done with much smaller hammers, but this is more complicated and might give the impression that we don't have big hammers. Also, as mentioned in the introduction, there is a rich history of partial results  not recounted here. 
\subsection{Resolutions with collared singularities}\label{ssect:collaredsing}
The proof given above uses the fact that singlarities in ANR homology manifolds can be ``resolved''. The next theorem gives a precise refinement for the PL case, based on the following lemma.
\begin{lemma}
Suppose \(L\) is a PL homology manifold  with the homology of a sphere. Then \(L\) bounds a contractible manifold in the sense that there is a contractible ANR homology manifold \(W\) with \(\partial W=L\), \(L\) has a collar neighborhood in \(W\), and \(W-L\) is a topological manifold. 
 Further, any two such \(W\) are homeomorphic rel a neighborhood of the boundary.
 \end{lemma}
 The only novelty is that we have not assumed \(L\) is a manifold. The proof of the 4-dimensional case given in \cite{fq}, Corollary 9.3C extends easily. We sketch the proof.
 
The standard triangulation of \(L\times[0,1]\) has no vertices in the interior so, by the Theorem above, the interior is a manifold. Do the plus construction (\cite{fq} \S11.1) to kill the fundamental group. The result is \(M\) with manifold interior, collared boundary \(L\times\{0,1\}\), and proper homotopy equivalent to a sphere. Replace each \(L\times[n,n+1] \subset L\times[0,\infty)\) by a copy of \(M\) and denote the result by \(W\). If \(W\) is a manifold except at the singularities of \(L\) then the standard manifold proof shows that the 1-point compactification is contractible, and a manifold except for these same singularities. It also shows that this manifold is unique up to homeomorphism rel boundary, The modification required in the older proof is verification that the interior of \(W\) is a manifold. 

 \(W\) is a manifold except possibly at vertices in  \(L\times\{n\}\) where the copies are glued together. If \(n>0\) then \(L\times\{n\}\) has a collar on each side, so has a neighborhood homeomorphic to \(L\times\R\), which is a manifold. Thus the only non-manifold points are in  \(L\times\{0\}\). This completes the proof.

 We use the Lemma to define models for ``collared singular points''. Suppose \(W\) is as in the lemma, with boundary collar \(L\times[0,1)\to W\). Identify the complement of a smaller open collar to a point to get \(W\to W/(W-L\times[0,1/2))\). The quotient is the cone \(L\times[0,1/2]/(L\times\{1/2\})\), the map is a homeomorphism except at the cone point, and the preimage of this point is a smaller copy of \(W\) and therefore contractible. In particular this is a resolution. 

Now define a ``resolution with collared singular points'' to be \(M\to K\) that is a homeomorphism except at a discrete set of points in \(K\), and near each of these points is equivalent to a standard model. The lemma easily implies:

\begin{theorem} A PL homology manifold \(K\) has a topological resolution with collared singular points, and  singular images the \(\pi_1\)-bad vertices of \(K\). This resolution is well-defined up to homeomorphism commuting with the maps to \(K\).
\end{theorem}
The mapping cylinder of a resolution is a homology manifold, and can be thought of as a ``concordance'' between domain and range. In these terms the theorem asserts that a PL homology manifold is concordant in a strong sense to a manifold. 

The unusually strong uniqueness (commuting exactly with maps to \(K\), not just arbitrarily close) results from the fact that two such resolutions have the same singular images, and the uniqueness in Lemma \ref{ssect:collaredsing}. This statement is true for manifolds with boundary if the definition of ``collared singularity'' is extended in the straightforward way.

\section{Triangulation}\label{sect:cssinvariant}
The main theorem is stated after the obstruction group is defined. The proof has two parts: first, enough structure of homology manifolds is developed to see the Cohen-Saito-Sullivan invariant. Both cohomology and homology versions  are described, in part to clarify the role of orientations. The second part is the converse, due to Galewski and Stern.

\subsection{The group}\label{ssect:theta3}

\(\Theta\) is usually defined as the set of oriented homology 3-spheres modulo homology \(H\)-cobordism. Connected sum defines an abelian monoid structure, and this is a group because reversing orientation gives additive inverses.   As mentioned in the introduction, the full name of this group is \(\Theta^H_3\), but analogous groups \(\Theta^H_k\) for \(k\geq 3\) are, fortunately, trivial. Roughly speaking, nontriviality would come from fundamental groups, and in higher dimensions we can kill these (eg.~with plus constructions). 

Geometric constructions give disjoint unions of homology spheres, not single spheres. These can be joined by connected sum to give an element in the usual definition of the group, but there are a number of advantages to using a definition that accepts disjoint unions directly. 
In this view \(\Theta\) is a quotient of the free abelian group generated by homology 3-spheres. Elements in the kernel are boundaries of oriented PL 4-manifolds that are homologically like \(D^4\) minus the interiors of finitely many disjoint 4-balls. These boundaries are disjoint unions of homology 3-spheres, and we identify disjoint unions with formal sums in the abelian group. 
Elements of the standard version are generators in the expanded version. It is an easy exercise to see that this inclusion gives an isomorphism of groups.  

In either definition it is important that the equivalence relation be  defined by PL manifolds, not just  homolgy manifolds. The goal is to organize \(\cone(L)\)-type singularities, and allowing singularities in the equivalences would defeat this. There may eventually be applications in which ``concordances'' can have singularities and the corresponding obstruction group should have these singularities factored out. For instance Gromov limits of Riemannian manifolds with special metrics might allow variation by cones on homology spheres with special metrics.

The Rokhlin invariant is a homomorphism \(\rok\colon \Theta\to Z/2\) defined using signatures of spin 4-manifolds bounding homology 3-spheres, cf.~\cite{kirby}. This connects with the Kirby-Siebenmann invariant, as described next.

\begin{theorem} (Main theorem)\label{thm:main}
\addcontentsline{toc}{subsection}{\ref{thm:main} \quad Main Theorem}
\begin{enumerate}\item \emph{(CSS invariant)} A PL homology manifold \(K\) has a `Cohen-Saito-Sullivan' invariant \(\css(K)\in H^4(K;\Theta)\);
\item \emph{(Relation to Kirby-Siebenmann)} If  \(r\colon M\to K\) is a topological resolution of a PL homology manifold then 
\(\ksm(M)=\rok(r^*(\css(K))\); and
\item \emph{(Realization: Galewski-Stern \cite{gs})} Suppose \(M\) is a topological manifold, not dimension 4, and a homeomorphism \(\partial M \to L\) to a polyhedron is given. If there is a lift \(\ell\) of  \(\ksm(M)\) to \(\Theta\) that extends \(\css(L)\) then there is a polyhedral pair \((K,L)\) and a homeomorphism \(M\to K\) that extends the homeomorphism on \(\partial M\), and \(\css(K)=\ell\).  
\end{enumerate}
\end{theorem}
 
Here, a ``lift'' is an element \(\ell\):
\[\xymatrix{ {\ell}\ar[r]\ar[d] &{\css(\partial K)}\ar[d]& \text{in}&H^3(M;\Theta)\ar[r]^{\partial^*}\ar[d]^{\rok}&H^3(\partial K;\Theta)\ar[d]^{\rok}\\
\ksm(M)\ar[r]&\ksm(\partial M)&&H^3(M;Z/2)\ar[r]^{\partial^*}&H^3(\partial K;Z/2)
}\]
There is a slightly sharper version in which ``lift'' is interpreted as a cochain representing such a cohomology class. Another extension is that if the map \(\partial M \to L\) is a resolution instead of a homeomorphism, then the conclusion is that it extends to \(M\to K\) that is a homeomorphism on \(M-\partial M\). The significance is that vertex singularities in \(L\) (where \(\partial M\to L\) cannot be a homeomorphism) do not effect the codimension-4 obstructions. Finally, the fact that 4-manifolds are smoothable in the complement of points \cite{ends3} can be used to alter definitions to give a formulation that includes dimension 4. We await guidance from applications to see which of these refinements is worth writing out.   

The proof of parts (1) and (2) are given in the remainder of this section. 
The Galewski-Stern proof of (3) follows the pattern developed to classify smooth and PL structures \cite{ks}, so is more elaborate than really needed. I did not find a proof short enough to include here, however.  
\subsection{Structure of polyhedra} We review the structure of polyhedra needed for homology manifolds. 
Suppose \(\sigma\) is a simplex in a simplicial complex. The \emph{dual cone} of \(\sigma\) is a subcomplex of the barycentric subdivision of the complex. Specifically, it is the collection of simplices that intersect \(\sigma\) in exactly the barycenter. The \emph{link} is the subcomplex of this consisting of faces opposite the barycenter point. 

\begin{figure}
\centering
\includegraphics[scale=0.6]{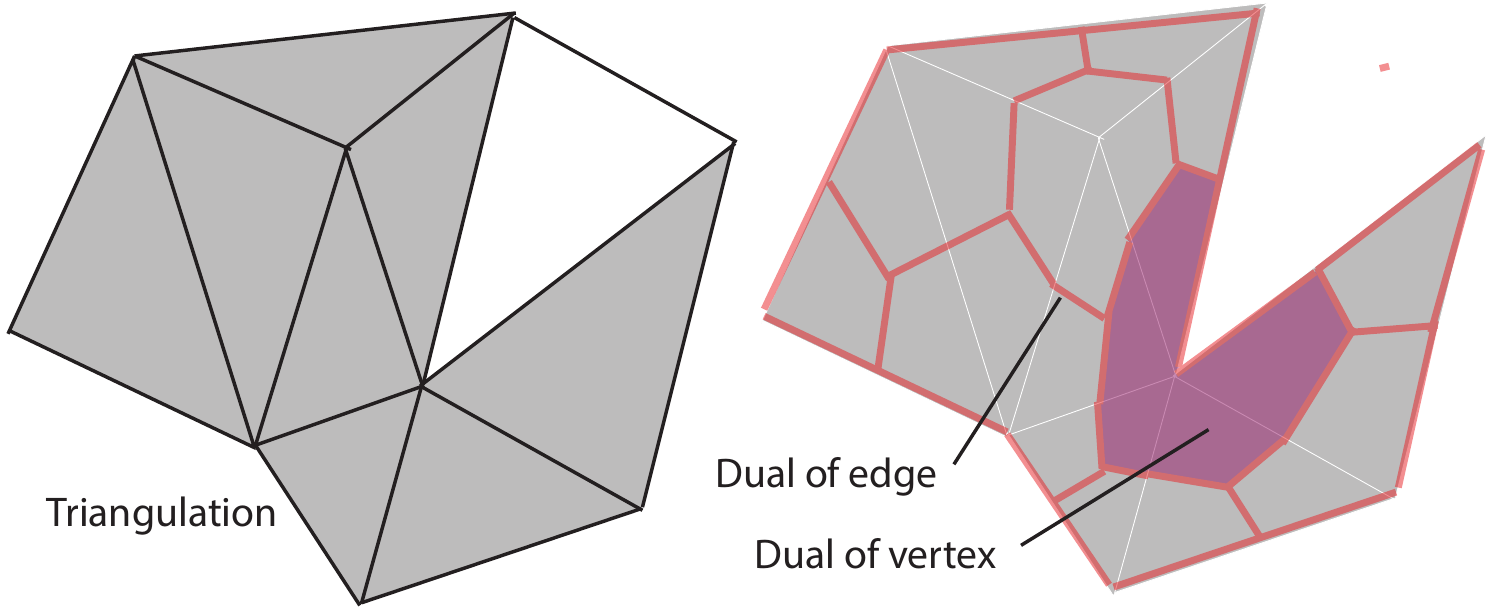} 
\caption{Dual cones in a simpicial complex}\label{fig:dualcells}
\end{figure}

It is easy to see that the dual cone is the cone on the link, with cone point the barycenter of \(\sigma\). The linear structures in the original simplices extends this to an embedding of the join of the link and \(\sigma\). Here we only need the weaker concusion that the interior of \(\sigma\) has a neighborhood isomorphic to the product \(\int (\sigma)\times \cone(\link(\sigma))\). 

\subsection{Links in PL homology manifolds}
Recall that  \(X\) is a homology \(n\)-manifold (without boundary) if  for each \(x\in X\), \(H_*(X,X-x;\Z)\simeq H_*(\R^n,\R^n-o;\Z)\). A pair  is a homology manifold with boundary if \(X-\partial X\) is a homology \(n\)-manifold, \(\partial X\) is a homology \((n-1)\)-manifold, and points in the boundary have the same local homology as points in the boundary of an \(n\)-ball (i.e.~trivial). 

PL homology manifolds have much more structure. 
\begin{lemma}\label{linklemma} A polyhedron \(K\) is a homology \(n\)-manifold (without boundary) if and only if link of every simplex is a homology manifold, and  has the homology of an \((n-k-1)\)-sphere, where \(k\) is the dimension of the simplex. 
\end{lemma}
The statement about homology of links is an easy suspension argument. The assertion that links are homology manifolds follows from this and the fact that links in a link  also appear as links in the whole space (easy after unwinding definitions).  This statement is easily extended to a version for manifolds with boundary. 

\subsection{The cohomology picture}\label{ssect:css}
 We begin with the cohomological version of the codimension-4 invariant. 

In a homology manifold the cones have the relative homology of disks, so they give a model for the chain complex. Specifically, define the \emph{conical chain} group \(C_n^{\cone}(K)\)  to be the free abelian group generated by \(n\)-dimensional cones together with a choice of orientation. 
Boundary homomorphisms in this complex come from homology exact sequences in a standard way.

Define a homomorphism \(Z[\text{oriented 4-d cones}] \to \Theta\) by
 \[(\cone(L),\alpha)\mapsto [L,\partial\alpha]\] where \(\alpha\) denotes the orientation of the (4-dimensional) cone, and \(\partial\alpha\) the corresponding orientation of the (3-dimensional) homology sphere. It is not hard to see that this defines a cohomology class \cite{cohen}, and we denote it by \(\css(K) \in H^4(K;\Theta)\). 
 
This definition includes manifolds with boundary, and the invariant of the boundary is \(\css(\partial K)=i^*\css(K)\), where \(i^*\colon H^4(K)\to H^4(\partial K)\) is induced by inclusion. The key result is that the Rokhlin homomorphism relates the Kirby-Siebenmann and CSS invariants:
\begin{proposition} If \(r\colon M\to K\) is a manifold resolution of a PL homology manifold then \(r^*(\rok(\css(K))=\ksm(M)\).
\end{proposition}
The usual formulation is for the special case with  \(r\) a homeomorphism. The additional information in the resolution version is that including the vertex singularities makes no difference. They neither contribute additional problems, nor do they give a way to avoid any of these problems.  

The description of \(\css\) should make this result very plausable, and if the definition of the Kirby-Siebenmann invariant is understood 
 (which we won't do here), the homeomorphism version should be obvious. The resolution version follows easily from the homeomorphism version and the  uniqueness of resolutions up to homeomorphism. 

\subsection{The homology picture}\label{ssect:HomPicture}
The dual homology class is sometimes easier to work with but takes more care to define correctly. The basic idea is to use simplicial chains, and represent the class in \(Z[(n-4)\text{-simplices}]\otimes \Theta\) by using the class of the link of a simplex \(\sigma\) as the coefficient on \(\sigma\). There is a problem with this: an orientation is required to define an element in \(\Theta\), but the data provides an orientation for the simplex rather than the dual cone. An orientation of the manifold can be used to transform simplex orientations to dual-cone orientations, but being too casual with this invites another mistake: the invariant is in twisted homology.

A homology manifold has a double cover with a canonical orientation, \(\hat K\to K\).  The group of covering transformations is \(Z/2\) and the generator acts on \(\hat K\) by interchanging sheets and (therefore) reversing orientation. Consider the simplicial chains \(C^{\Delta}(\hat K)\) as a free complex over the group ring \(Z[Z/2]\), and suppose \(A\) is a \(Z[Z/2]\) module. We define the homology \(H_n(\hat K; A)\) to be the homology of the complex \(C^{\Delta}(\hat K)\otimes A\), where the tensor product is taken  over \(Z[Z/2]\). 

If \(Z/2\) acts trivially on \(A\) then the tensor product kills the action on the chains of \(\hat K\) and the result is ordinary homology. We will be concerned with the opposite extreme, \(A=Z\) with \(Z/2\) acting by multiplication by \(-1\). 

After this preparation we can define the Cohen-Saito-Sullivan \emph{homology} class  by
\[\css_*({K}) = \Sigma_{\sigma} \sigma\cdot[\link(\sigma)] \in H_{n-4}(\hat K, \partial \hat K; \Theta)\]
where \(Z/2\) acts on \(\Theta\) by reversing orientation, and the orientation of \(\link(\sigma)\) is induced by the orientation of \(\sigma\) and the canonical orientation of \(\hat K\). 

This definition also includes manifolds with boundary, and the invariant of the boundary is given by the boundary homomorphism in the long exact sequence of the pair.

The homology and cohomology definitions are Poincar\'e dual. Duality between simplices and dual cones is particularly clear: each simplex intersects exactly one dual cone (its own) in a single point, and this pairing gives a chain isomorphism between simplicial homology and dual-cone cohomology when links are homology spheres\footnote{This is, in fact, Poincar\'e's picture of duality, and will be discussed further in the history section.}. This pairing matches up the two definitions. 

\section{Gauge theory}\label{sect:gaugetheory}
The Casson invariant (see \cite{am}) gave first hint that something like gauge theory would play a role in this story. Casson used representation varieties and Heegard decompositions to define an integer-valued invariant of homology 3-spheres, and showed that the mod 2 reduction is the Rokhlin invariant. However it is a invariant  of  diffeomorphism type,  not homology H-cobordism. It does not define a function \(\Theta\to\Z\), and has little consequence for the triangulation questions. 

Fintushel and Stern \cite{fintstern90} used the Floer theory associated to Donaldson's anti-self-dual Yang-Mills theory to show that certain families of Seifert fibered homology 3-spheres are linearly independent in \(\Theta\). The families are infinite so \(\Theta\) has infinite rank. This implies that most manifolds have vastly many concordance classes of triangulations, but does not clarify the existence question because all  these homology 3-spheres are in the kernel of the Rokhlin homomorphism. 

There has been quite a bit of work done since Fintushel-Stern, with invariants derived from gradings in various Floer homology theories; see  Manolescu's discussion of Fr\o yshov correction terms. The next qualitatively new  progress, however, is in Manolescu's paper. The outcome is three functions \(\Theta\to\Z\) which are not homomorphisms, but have enough structure to show that a homology sphere with nontrivial Rokhlin cannot have order 2 in \(\Theta\). This implies that there are manifolds of dimension 5 and higher that cannot be triangulated: those whose Kirby-Siebenmann mod 2 classes do not lift to mod 4 cohomology. See \cite{gs5} for a 5-dimensional example. Somewhat more elaborate arguments with these functions seem to show that many Rokhlin-nontrivial spheres have infinite order. The full consequences are not yet known. 

Sections \ref{ssect:physics}--\ref{ssect:equiv} gives a qualitative outline of Manolescu's paper. References such as ``\cite{mano13} \S3.1'' are abbreviated to `M3.1', and readers who want to see things like the Chern-Simons-Dirac functional written out should refer to this paper. Alternate perspectives for experts are suggested in  \S\ref{ssect:handcraft}--\ref{ssect:goals}

\subsection{Physics description}\label{ssect:physics}
The idea, on the physics level of rigor, is that the Floer homology theory associated to the Seiberg-Witten equations is given by the Chern-Simons-Dirac functional on an appropriate function space. This functional is invariant under a big symmetry group. Divide by the symmetry group, then we want to think of the induced function on the quotient as a sort of Morse function and study its gradient flow. More specifically, we are concerned with the finite-energy trajectories.  The quotient is infinite-dimensional, but we can enclose the finite-energy trajectories in an essentially finite-dimensional box. Invariants of the system come from algebraic-topological invariants of this box. 

This description  offers an alarmingly large number of ways to misunderstand the construction, and one goal is to clarify the strategy and logical structure of the process. For instance, finite-dimensional differential and algebraic topology are mature subjects with a lot of sharp tools. It is useful to see the infinite-dimensional part of the analysis as a sequence of reductions designed to bring part of the structure within range of these sharp tools.   
Another, possibly dubious, goal is to try to clarify features of the technical details and how they might be sharpened, but without actually describing the details. Finally, the analysis described is for  3-manifolds whose  first homology  is  torsion (\(b_1=0)\). The analysis in the general case is considerably more elaborate. 
\subsection{The Coulomb slice}\label{ssect:coulomb}
The first step in the heuristic description is to  ``divide by an infinite group of symmetries''. It is almost impossible to make literal sense of this, and in    M3.1 Manolescu uses the Coulomb slice to avoid it.  There is a (``normalized'') subgroup of the full symmetry group with the property that each orbit intersects this slice in exactly one point. The slice is therefore a model for the quotient by this subgroup, and projection to the slice reduces the symmetry to the quotient of whole group by the subgroup.  The quotient is the compact Lie group \(Pin(2)\). 

Since it is compact, dividing by   \(Pin(2)\)  makes good sense, but it introduces singularities that are much more painful than symmetry groups.  The plan is therefore to do a nonsingular equivariant reduction to finite dimensions, and the long-term strategy is roughly ``let the finite-dimensional people deal with the  group action''. 

Manolescu explicitly describes the restriction of the Chern-Simons-Dirac functional to the Coulomb  slice, and describes a projected Riemannian metric that converts the derivative of the functional to a gradient vectorfield with the property that the projection preserves gradient flows. This description usually gives non-specialists the wrong picture because the ``Riemannian metric''  is not complete. The slice is a Frech\'et space of \(C^{\infty}\) functions, and there is no existence theorem for flows in this context. In fact, in most directions the gradient vectorfield does not have a flow, even for short time, and the flow trajectories exploited by Floer and others exist due to a regularity theorem for solutions of a differential equation with boundary conditions. In other words there is only a small and precious fragment of a flow for this vectorfield, and  this is \emph{not}  Morse theory with a globally-defined flow. The observation that projection to the Coulomb slice preserves flows means it preserves this small and precious fragment, not something global. 

This explanation is still not quite right. Manolescu doesn't actually identify the flow fragment in infinite dimensions, so saying that the projection preserves whatever part of the flow that happens to exist is a heuristic summary. On a technical level the projection preserves \emph{reasons} the fragment exists and it is these reasons, not the flow itself, that power the rest of the argument.

\subsection{Sobolev completions}\label{ssect:sobolev} In the last paragraph of M3.1 the space \(V_{(k)}\) is defined as the completion of the Coulomb slice, using the \(L^2\) Sobolev norm on the first \(k\) derivatives. This gives Hilbert spaces but still doesn't give us a flow because the ``vectorfield'' now changes spaces: it is of the form 
\[\ell +c\colon V_{(k+1)}\to V_{(k)}\]
with \(\ell\) linear Fredholm and \(c\) compact\footnote{The spaces \(V\) are vector spaces and \(\ell\) linear because we are assuming \(b_1=0\) (homology sphere). The general situation is more complicated.}.  The index shift corresponds to a loss of a derivative, reflecting the fact that we are working with a differential equation. 
The manouvering (bootstrapping) needed to more-or-less recover this lost derivative is a crucial analytic ingredient. Almost nothing is said about this in \cite{mano13}, but some details are in \cite{mano03}, sections 3 and 4, phrased in terms of flows rather than vectorfields.   Manolescu's next step is projection to finite dimensional spaces where there are well-behaved flows. There would be significant advantages to connecting directly with Morse theory in an infinite-dimensional setting rather than in projections; see  \S\ref{sssect:hilbert} for further comments.  

\subsection{Eigenspace projections}\label{ssect:proj} 
In section M3.2 Manolescu defines \(V^{\nu}_{\tau}\) to be the subspace of \(V\) spanned by eigenvectors of  \(\ell\) with eigenvalues in the interval \((\tau, \nu]\). This uses the fact, prominent in \cite{mano03} but unmentioned in \cite{mano13}, that \(\ell\) is  self-adjoint.  In particular its eigenvalues are real and eigenspaces are spanned by eigenvectors. These spaces are finite dimensional because \(\ell\) is Fredholm. Finally, the symmetry group \(Pin(2)\) acts on them because they are defined using equivariant data.

There is a technical modification that deserves comment. The orthogonal projections \(V\to V^{\nu}_{\tau}\) give a function from the parameter space \(\{\tau<\nu\}\) to linear maps \(V\to V\). This takes discrete values (depending only on the eigenvalues in the interval) so is highly discontinuous. Manolescu smooths this function: the dimension of the image still jumps but when it does, the projection on the new part is multiplied by a very small number. The result is continuous as a function into the space of linear maps. This implies that the finite projections of the CSD vectorfield become smooth functions of the eigenvalue parameters. This is useful in showing that parts of the qualitative structure of the output do not depend on the parameters once they  are sufficiently large.

 The final modification of the flow is done in section M3.7. There is a unique reducible solution of the equations, and non-free points of the \(Pin(2)\) action come from this. The functional  is perturbed slightly (following the earlier  \cite{mano03}) to make the reducible solution a nondegenerate critical point. The irreducible critical points can also be made nondegenerate in an appropriate equivariant (Bott) sense. 

This is one of the places where the Seiberg-Witten theory diverges in a qualitative way from the Donaldson theory. The finite-energy trajectories in the Donaldson-Floer theory cannot be made nondegenerate, and the analysis takes place on a center manifold. This is rather more delicate. 

\subsection{Isolated invariant sets} The last structural input from the infinite dimensional context is specification of the ``precious fragment'' of the flow supposed to have come from infinite dimensions. This is done by Proposition M3.1 in \cite{mano13}, which is a reference back to Proposition 3 of \cite{mano03}. The flow fragment is the union of trajectories that stay in a ball of a certain radius, and the key fact is that it is isolated in the sense that it is the same as the union of trajectories that stay in a ball of twice the radius. 

We comment on the logic of the reduction. Defining the invariant uses only the answer (the form of the explicit finite dimensional approximations) and the proof in Proposition M3.1 that the trajectories-in-a-ball construction gives an isolated invariant set. This does not use the construction of a flow fragment in infinite dimensions, so the demonstration that such a flow fragment would have been preserved by the projection is not actually used. This demonstration does, however, give a tight connection between this construction and those of Floer et.al.~that do use the infinite-dimensional flow. 

\subsection{Equivariant stable homotopy theory}\label{ssect:equiv}  The plan is to enclose the isolated invariant set identified in the previous step in a nice box, and extract information about the system from algebraic and geometric topology of the box. The box is a subspace (or submanifold) of a finite-dimensional vector space so this is the point at which the problem enters the finite-dimensional world. Manolescu is not a native of this world, however, and his treatment could be refined. We briefly sketch  Manolescu's definition of the invariants in this section. The main difficulties come in showing that these are well-defined and have good properties. The next section hints at some of these difficulties and suggests approaches that may be better adapted. 

Manolescu uses the Conley Index construction to get a ``box'' enclosing the isolated invariant set.  
The output is a pair of spaces that depends on choices, or by taking the quotient a pointed space that the choices change only by homotopy equivalence. This takes place in an eigenvalue projection   \(V^{\nu}_{\tau}\) and changing \(\tau\) and \(\nu\) changes the pointed-space output by suspension. The object associated to the system is therefore a spectrum in the homotopy-theory sense. Finally, all these things have  \(Pin(2)\) actions, and the suspensions are by \(Pin(2)\) representations.  The proper setting for all this is evidently some sort of equivariant stable homotopy theory. The most coherent account in the literature is  Lewis-May-Steinberger \cite{may86}, and Manolescu uses this version. 
The next step is to extract numerial invariants from these  \(Pin(2)\)-equivariant spectra using Borel homology.

To a first approximation the homology appropriate to a \(G\)-space \(X\) is the homology of the quotient \(X/G\). This works well for free actions but undervalues fixed sets. The Borel remedy is to make the action free by product with a contractible free \(G\)-space \(EG\), and take the homology of the quotient \((X\times EG)/G\). The free part of \(X\) is unchanged by this but  points fixed by a subgroup \(H\subset G\) are blown up to copies of the classifying space \(EG/H\). These classifying spaces are usually homologically infinite dimensional, so fixed sets become quite prominent. Another benefit of the Borel construction is that the homology of  \((X\times EG)/G\) is a module over the cohomology of \(BG\colon=EG/G\).  These facts are illustrated by a localization theorem quoted in M2.1\footnote{The finiteness hypothesis on \(X\) is missing in the statement in \cite{mano13}.}: suppose \(X\) is a finite \(G\)-complex and the action is free on the complement of \(A\subset X\). A localization that kills finite-dimensional \(H^*(BG)\) modules  kills the relative Borel homology of \((X,A)\),  so the inclusion  \(A\to X\) induces an isomorphism on localizations. There is a difficulty that Borel homology is not fully invariant under equivariant suspensions. Manolescu finesses this with \(F_2\) coefficients, but eventually it must be addressed. 

In the case at hand the \(G\)-objects are spectra rather than single spaces. \(X\) can be thought of as the equivariant suspension spectrum of a finite \(G\)-complex and the sub-spectrum \(A\) of non-free points is essentially is the suspension spectrum of a point. Inclusion therefore gives a \(H^*(BG)\)-module homomorphism  \(H_*(BG)\to H_*((X\times EG)/G)\). Finiteness of \(X\) implies that the third term in the long exact sequence (the homology of the free pair \((X,A)\)) is finite-dimensional. In particular the kernel of \(H_*(BG)\to H_*((X\times EG)/G)\) is a finite-dimensional \(H^*(BG)\)-module.  When \(G=Pin(2)\) these submodules are characterized by three integers \(\alpha, \beta, \gamma\), and these are Manolescu's invariants. The algebraic details give a pretty picture, and readers should refer to Manolescu's paper for this.

\subsection{Handcrafted contexts}\label{ssect:handcraft} Both stable homotopy theory and equivariant topology are sprawling, complicated subjects. Off-the-shelf versions tend to be optimized for particular applications and often use shortcuts or sloppy constructions that can cause trouble in other circumstances.  The best practice is to handcraft a theory that fits the application, but this requires insider expertise. In this section we suggest such a handcrafted context for the finite-dimensional part of Manolescu's development. 
\subsubsection{Lyapunov blocks}
The first step is to be more precise about the data at the transition from analysis to finite-dimensional topology. Manolescu uses the Conley index construction to get a ``box'' enclosing an isolated invariant set in a flow on a manifold \(V\). We recommend instead an object we call a \emph{Lyapunov block}. These were introduced and shown to exist using Lyapunov functions by Wilson and Yorke \cite{wilsonyorke73}, and shown to be essentially equivalent to Lyapunov functions by Wilson \cite{wilson80}. Wilson and Yorke call these ``isolating blocks'', but a more distinctive name seems to be needed. This construction has been revisited recently by Cornea \cite{cornea}, Rot-Vandervorst \cite{rotvdh}, and others.

\begin{figure}
\centering
\includegraphics[scale=0.7]{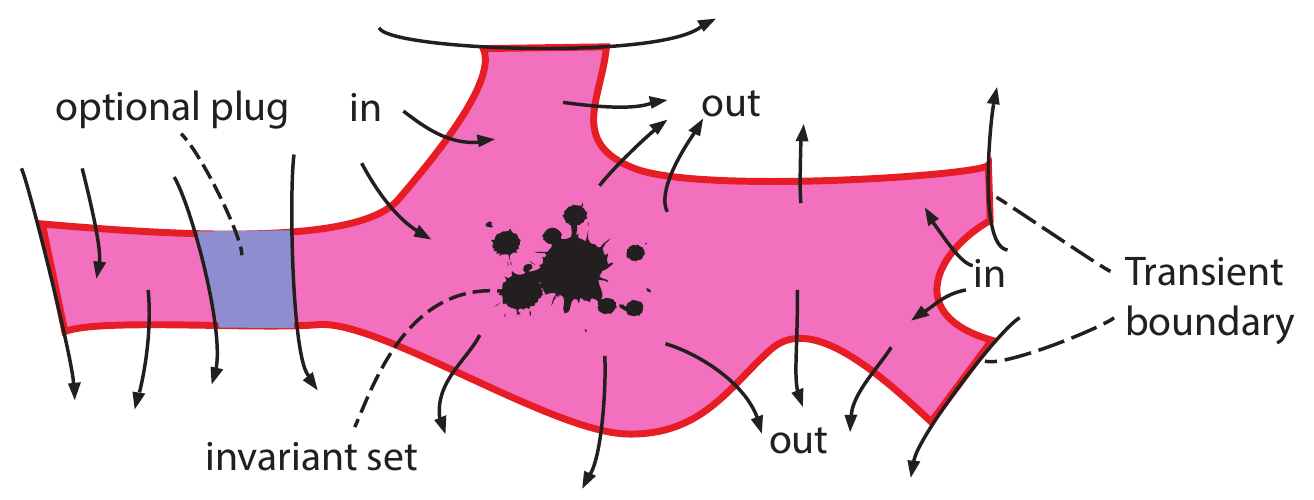} 
\caption{A Lyapunov block for a flow}\label{fig:lyapunovbox}
\end{figure}

A  Lyapunov block for an invariant set in a flow is a compact smooth codimension-0 submanifold with corners  \(B\subset V\) with boundary divided into submanifolds \(\partial_-B\cup \partial_0B\cup  \partial_+B\). Trajectories intersect \(B\) in arcs.  Trajectories enter through the incoming boundary \(\partial_+B\), exit through the outgoing boundary \(\partial_-B\). The transient boundary \(\partial_0B\) is a union of intersections with trajectories, and the trajectory arcs give a product structure \(\partial_0B\simeq\partial_{0,+}B \times I\); see Figure \ref{fig:lyapunovbox}. Finally, the trajectories completely contained in \(B\) are those in the original invariant set.  The underlying smooth manifold structure  can be thought of as a smooth manifold triad \((B,\partial_+B,\partial_-B)\).   \(\partial_0B\) is a collar so absorbing it into either \(\partial_-B\) or \(\partial_+B\) (or half into each) changes  them only by canonical diffeomorphism. 

These blocks are not well-defined: different choices in a truncation step give \(B\) that differ by addition or deletion of plugs of the form \(P\times I\), that intersects trajectories  in product arcs \(\{p\}\times I\). This implies that the pairs \((B,\partial_+B)\) and \((B,\partial_-B)\) have well-defined relative homotopy types\footnote{Blocks can be modified to eliminate the transient boundary \cite{rotvdh}, but it is best not to make this part of the definition because it makes the ``plug'' variation hard to formulate.}.  
To relate this to Manolescu's version,  \((B,\partial_-B)\) is a particularly nice Conley index pair for the flow, and the index itself is the pointed space \(B/\partial_-B\). The quotient \(B/\partial_+B\) is a Conley index for the reversed flow. The manifold triad  therefore gives both Conley indices and precisely encodes their relationship.
\subsubsection{Smooth manifold triads}
The handcrafted context appropriate to this situation seems to be a stable category defined using equivariant smooth manifold triads. It is ``stable'' in the sense that the objects are families of triads related by equivariant suspensions that are ``internal'' in the sense that they come from eigenvalue-range changes (see below).  This context receives  Lyapunov blocks  without further processing. It has many other virtues, as we explain next, and in fact we like Lyapunov blocks because they permit use of this context. 

This context does \emph{not} follow the standard practice of dividing to get pointed spaces. Data from geometric situations often comes as pairs with structure that  does not gracefully extend to pointed-space quotients. Bundles on pairs, for instance, rarely extend to the pointed space. This means they have to be described as ``bundles over the complement of the basepoint'', and to work with them one must recover the pair by deleting a neighborhood of the basepoint. Group actions can be extended to have the quotient basepoint as a fixed point, but this is often just cosmetic. In many geometric applications, for instance, algebraic topology is done equivariantly on the universal cover. The fact that the action is free is essential. The pointed-space quotient therefore must be described as an action free in the complement of the basepoint, and again much of the work requires deleting the basepoint to recover a pair with a free action. Having to delete the basepoint is a clue that dividing to get a basepoint was a mistake. In some cases the pair information can be recovered stably without expicitly deleting the basepoint, but it is usually a lot of work.

The second advantage of this context is that in the manifold-triad world, Spanier-Whitehead duality is implemented by interchanging the two boundary components. Interchanging boundaries in a Lyapunov block corresponds to reversing the flow, so it is obvious that the flow and its reverse have S-W dual blocks.  In the pointed-space context there is a stable description in terms of maps from a smash product to a sphere, but this is a characterization, not the definition, and it does not work in all cases. Manolescu quotes this version in M2.4, and cites references that show that the stable and unstable Conley indices are S-W dual in this sense. But these references use Lyapunov blocks, so the net effect is ``discard the Conley constructions and redo the whole thing with manifold triads''. Going back to Conley indices not only is inefficient but introduces troublesome ambiguities about suspensions. This difficulty is discussed next. 

The final wrinkle in this context has to do with the meaning of ``suspension''. Enlarging the range of eigenvalues changes the projection by product with a representation of \(Pin(2)\), and changes the Lyapunov block  by suspension with the ball in this representation. This is an ``internal'' suspension because it is specified by the analytic data. Understanding how internal suspensions change, for instance when the metric on the original manifold is varied,  is a job for analysis. External suspensions used to define equivariant invariants are specified differently, and the two types of suspensions should be kept separate. In particular the eigenvalue-change suspensions should not be seen as instances of external suspension operations. To explain this, note that the equivariant theory of Lewis, May and McClure \cite{may86} (used by Manolescu) is handcrafted to give a setting for homology theories and classifying spaces. Roughly speaking, they want to grade homology theories by equivalence classes of objects in the category of representations. When objects have nontrivial automorphisms, equivalence classes of objects do not form objects in a useful category. The standard fix for this is to use a skeleton subcategory  with one object in each equivalence class. In the equivariant setting this means choosing one representation in each equivalence class, and always suspending by exactly \emph{this} representation. This is fine for \emph{external} suspensions, but 
representations that come internally from eigenvalue projections have no canonical way to be identified with randomly chosen representatives. If the group is \(S^1\), as in most previous work on Seiberg-Witten-Floer theory, then there are essentially no automorphisms and this issue can be finessed. Manolescu's key insight, however, is that \(Pin(2)\) is the right symmetry for this problem\footnote{The \(Pin(2)\) symmetry was observed much earlier, cf.~\cite{bauerFuruta}, but not fully exploited.}, and these representations have automorphisms that make  identifications problematic. The solution is to avoid using external suspensions in describing the geometric invariant. Lyapunov blocks do this.

\subsection{Next goals} \label{ssect:goals}
Floer homology is, to a degree, a solution in search of worthy problems. Distinguishing knots is a baby problem whose persistance just reflects the lack of real work to do. The triangulation problem is useful for teething  technology but, as explained in the History section, is a backwater with no important application. This weakness is reflected in the structure of Manolescu's invariant: detailed information about homology spheres lies in the part of the moduli space on which \(Pin(2)\) acts freely, but the invariant discards all this except the levels at which it cancels homology coming from the fixed point. This is not a gateway to something deeper. We have several suggestions for further work. 

\subsubsection{Complicated geometric structure}
The first suggestion is motivated by  internal structure of the analytic arguments. Analysis  associates to a  homology 3-sphere  a complicated \(Pin(2)\)-equivariant gadget. We expect this to reveal something about geometric properties of the 3-manifold, but neither the properties nor the mechanisms of revelation are clear. A useful intermediate step would be a complicated \(Pin(2)\)-equivariant gadget derived more directly from the topological object. The two equivariant gadgets might be related by a sort of index theorem. The point is that sometimes it is easier to relate two complicated things than to understand either in detail, and the connection can be a powerful aid to understanding.  We suggest maps from the homology sphere to \(S^3\simeq SU(2)\) as the topological gadget. \(Pin(2)\) acts on this because it is a subgroup of \(SU(2)\), and these maps should connect to geometric structure by a form of generalized Morse theory, cf.~\cite{gaykirby}. 
\subsubsection{Hilbert, or SC manifolds}\label{sssect:hilbert}
The key analytic goal is to situate the objects of interest in a context accessible to  ``finite-dimensional'' geometric and algebraic topology. The context does not have to be literally finite-dimensional to use the techniques, however, and a context that does not require  finite-dimensional projections would simplify formulation of invariants. The first requirement for such a context is an effective global existence theorem for flows. There seem to be at  two  possibilities that are, in a sense,  at opposite extremes. 

Manolescu begins (see  \S\ref{ssect:coulomb}) by restricting the Chern-Simons-Dirac functional  to the Coulomb slice, and using a Riemannian metric to convert the derivative of the functional to a vector field. It would be quite natural to complete with respect to this metric, to get a vectorfield on a separable Hilbert manifold. The problem is that current estimates are not good enough to show that the finite-energy trajectories form an isolated invariant set in this topology. There are  heuristic reasons to worry that they are not isolated in general. A perturbation of the system to be `nondegenerate' in some sense might help.
 Eventually the geometric invariant would be a Lyapunov block in the Hilbert manifold, together with an equivalence class of \emph{structures} related to the finite-dimensional projections. This would clarify that the objects obtained by projection are fragments of a structure on the invariant object, not the invariant object itself. A Hilbert-manifold formulation should be much easier to extend to things like Hilbert-manifold bundles over \(H_1(X;R)\), which may be necessary for 3-manifolds with \(\beta_1\neq0\). 

 Another possible context is the Banach-scale manifolds  developed by Hofer, Wysocki, and Zehnder.  The Hilbert approach takes place at a fixed level of differentiability, while the Banach-scale approach organizes the way in which function spaces of increasing differentiability approach \(C^{\infty}\).  It is `handcrafted' in the sense of \S\ref{ssect:handcraft} to formalize and exploit the bootstrapping common in applications.
 
 The first comment is that the\ Hofer-Wysocki-Zehnder ``polyfold'' theory is not appropriate here. This was developed to handle closure problems in quotients. Here this is handled by taking the quotient by a subgroup of the full gauge group which, since it has a global slice, has no closure problems. This leaves a residual gauge action by \(Pin(2)\). Dividing by this does introduce orbit-closure problems but (1) these seem to be outside the reach of the Hofer-Wysocki-Zehnder polyfold theory, and (2) by now it should be quite clear that equivariant nonsingular objects are more effective  than trying to describe some sort of structure on singular quotients. The second comment is that a useful version of ``Lyapunov block'' would be needed, and this may require negotiation with the topological theory that has to use it. The final comment is that this may give a setting for the germ-near-a-compact-set suggestion in the next section. 
 
\subsubsection{Stay in dimension 4}
The motivation for the final suggestion is external to the analysis. The  best guides to development of a theory are deep potential applications.  Floer homology of 3-manifolds is supposed to organize boundary values and glueing properties of gauge theories on smooth 4-manifolds but, in general, 3-manifolds slices and boundaries do not adequately reflect the complexity of smooth 4-manifolds. We explain this in a context that ideally would connect with homology 3-spheres. 

Suppose \(M\) is a smooth 4-manifold with a submanifold \(V\) \emph{homeomorphic} to \(S^3\times \R\).  If \(M\) is compact, simply-connected, and \(V\) separates \(M\) then  a relatively soft argument \cite{freedmanTaylor} shows that \(M\) also contains a smooth homology 3-sphere homologous to \(S^3\times \{0\}\). But this is usually not true if \(M\) is either noncompact or not simply connected. For instance, a compact 4-manifold has a smooth structure in the complement of a point, and this point has a neighborhood homeomorphic to \(S^3\times \R\), but almost none of these contain  smooth homology spheres. When there is an appropriate homology sphere in \(M\) it is usually not in the given \(V\).

Another soft argument shows that in the compact simply-connected case any two  homology 3-spheres arising as above are  homology \(H\)-cobordant, but not ``in \(M\)''. Note that disjoint homologous homology spheres have a region between them that is an \(H\)-cobordism.   Ideally, if we have two homology spheres then  we would find  a third homology sphere disjoint from both. The first two would both be  \(H\)-cobordism to the third, so the first two would be \(H\)-cobordant by a composition of \emph{embedded} \(H\)-cobordisms. Unfortunately we can usually not find a third disjoint sphere, and the soft argument does not give embedded \(H\)-cobordisms. 
 
 The moral of this story is that we can use  transversality to get smooth 3-manifold splittings, but these 3-manifolds usually cannot reflect the global homotopy theory of the manifold even up to homology. A glueing theory that depends on finding nice slices (eg.~smooth homology spheres in topological connected sums)   therefore cannot be an effective general setting. 

A better setting for glueing 4-d theories should be some sort of ``germs of necks'' that locally separate the 4-manifold. We have much better criteria for finding good \emph{topological} slices in 4-manifolds, so a first approximation would be  ``germs near \(X\times \{0\}\) of smooth structures on \(X\times \R\)'', where \(X\) is a closed 3-manifold, but the smooth structure on  \(X\times \R\) is not the product structure.  

Smooth neighborhoods of topological embeddings is the sort of mixed-category thing that (according to the History section) is probably a bad idea in the long term, but it gives a concise starting point.  The homotopy data required to find a topological slice in a ``neck'' are non-obvious and fairly elaborate. The data needed to find a ``virtual analytic slice'' may also be elaborate, so speculations should wait on feedback from analysis. In any case the point for the present discussion is that the best next step in Floer-type theory is probably gauge theory on 4-d ``neck germs'', not gauge theory on 3-manifolds.

\section{History}\label{sect:history} 
Poincar\'e's  insights about the homology of manifolds, at the end of the nineteenth century, are usually celebrated as the starting point of modern topology. But many of his insights were wrong in detail, and his methodology was so deficient that it could not be used as a foundation for further development. His contemporaries found it inconceivable that the Emperor might have no clothes, so they spent the next quarter-century trying to see them. 
 Kneser's triangulation questions are precise formulations of what it would take to make Poincar\'e's  arguments sensible. Labeling one of them `the Hauptvermutung' suggests that he still hoped it would all work out. But it did not. As interesting as these questions seem, they are a technical dead end: not only not a foundation for manifold theory, but apparently without any significant applications. Details of this story, and how topology finally recovered from Poincar\'e's  influence, are told in the rest of this section. 

\subsection{Pre-modern methodology}
Poincar\'e worked during the period when modern infinite-precision mathematics was being developed \cite{rev}. He was not part of this development, however, but worked in---and strongly defended---the older heuristic and intuitive style. His explanations often included the technical keys needed for a modern proof of a modern interpretation of his assertion. But he often omitted hypotheses necessary for his assertions to be correct, and his arguments were too casual to reveal the need for these hypotheses. He gave examples, but did not use precise definitions and often did not verify that the examples satisfied the properties he ascribed to them. This casual approach, and the philosophical convictions that underlay it, made for a difficult start for the subject. 

 For instance Poincar\'e proceeded on the presumption that the choice of analytic, combinatorial, or topological tools would be dictated by the task at hand rather than the type of object. Functionally this amounts to an implicit claim that topological, PL, and smooth manifolds are all the same. Clearly anything built on this foundation was doomed. But identifying this as a flaw in Poincar\'e's work would have invited strong political and philosophical attack and the new methodologies were not secure enough for this.    Kneser's triangulation questions twenty five years 
 were  precise technical formulations of what would be needed to justify Poincar\'e's work and approach, but he still did not identify this as a gap in the work. 

Not only was it hard to know which parts of Poincar\'e's work were solid, but apparently it was hard to track which parts were actually known to be false.  For instance in a 1912 paper of Veblen and Alexander \cite{veblenAlex12} we find 
\begin{quote}
Poincar\'e has proved that any manifold \(M_n\) may be completely characterized from a topological point of view by means of suitably chosen matrices \dots
\end{quote}
This refers to the 1895 claim that the incidence matrices of a triangulation (now called boundary homomorphisms in the chain complex) characterized manifolds up to homeomorphism. We overlook this blunder today because Poincar\'e himself disproved it  not long after, by using the fundamental group  to show the ``Poincar\'e sphere'' is not \(S^3\) even though it has equivalent chains. But more than ten years later Veblen and Alexander seem to have been unaware of this refutation. 
\subsection{Poincar\'e's duality}
An explicit example of Poincar\'e's methodology is provided by his description of duality. He observed the beautiful pairing of simplices and dual cones explained in \S\ref{ssect:HomPicture}. But he called these dual cones `cells', and implicitly presumed that they were equivalent (in an unspecified sense) to disks. Instead of seeing this as a general PL construction that might or might not give a cell, it was seen as a manifold construction that ``failed'' if the output was not a cell.  This convention makes arguments with dual cells logically sensible, but it hides the necessity of showing that the construction does not fail in specific instances. One  of Poincar\'e's classes of examples was inverse images of regular values of smooth maps \(R^n\to R^k\). In what sense can these be triangulated, and why are the dual objects cells? Whitehead sorted this out some 40 years later \cite{whitehead40}. The proof was probably beyond Poincar\'e's ability, but the real problem was that he did not notice (or acknowledge) that there was a gap.

Another difficulty is that Poincar\'e's duality relates different objects: homology (or Betti numbers) based on simplices on the one hand, and homology based on dual cells on the other. In order to get duality as a symmetry of a single object, these must be identified in some other way.  This time Poincar\'e got it wrong: dual cells actually give cohomology so, as we know now, duality gives an isomorphism between homology and cohomology. The homology/cohomology distinction (in the group formulation) together with the Universal Coefficient Theorem explain why the torsion has symmetry shifted one dimension from the Betti-number symmetry. Poincar\'e missed this, and found a patch only after Heegard pointed out a contradiction. In another direction, duality requires some sort of orientation and (as we saw with \(\Theta\) in \S\ref{ssect:HomPicture}) may be twisted even when there is an orientation. When the manifold has boundary, or is not compact, duality pairs homology with rel-boundary or compact-support homology. Homology of the boundary appears as an error term for full symmetry.  Again these results were beyond Poincar\'e's  intuitive definitions and heuristic arguments, but the real problem was that he did not notice (or acknowledge) that more argument was needed. 

 \subsection{Point-set topology}  Schoenfliss and others were developing point-set topology around the same time, and the relationship between the two efforts is instructive.
 
An important point-set goal  was to settle the status of the Jordan Curve theorem. This  is not hard to prove for smooth or PL curves, but an intuitive extrapolation to continuous curves was discredited by the discovery of continuous space-filling curves by Peano and others. The continuous version had important implications for the emerging role of topology as a setting for analysis. For instance, integration along a closed curve around a ``hole'' in the plane was a vital tool in complex analysis. Integration required  piecewise-smooth curves. The question was: were ``analytic holes'' identified by piecewise smooth curves the same as ``topological holes''  identified with continuous  curves? If not then the role of general topology would probably be quite limited.
 
Addressing the Jordan Curve problem turned out to be  difficult, and fixing gaps in attempted proofs required  quite a bit of precision about open sets, topologies, separation properties, etc. In short, it required modern infinite-precision techniques. Wilder \cite{hist} found it curious that Shoenfliss never mentioned Poincar\'e or his work, since nowdays the Jordan Curve theorem and high-dimensional analogues are seen as immediate consequences of a homological duality theorem. But this makes  sense: Shoenfliss was trying to fix a problem in a heuristic argument, and  Poincar\'e used, and strongly defended the use of, heuristic arguments. The duality approach was not available to Shoenfliss because---for good reason---he could not trust Poincar\'e's statements about duality.  
\subsection{Constrained by philosophy} The general question in this section is: why did it take Poincar\'e's successors so long to find their way past his confusions? The short answer is that they were in the very early part of the modern period and still vulnerable to old  and counterproductive convictions. We go through some of the details for what they reveal about the short answer: what \emph{were} the counterproductive nineteenth-century convictions, and how did they inhibit mathematical development?

To be more specific, a  mathematician with modern training would probably respond to  Poincar\'e's work  with something like\begin{quotation} The setting seems to be polyhedra, and the key property seems to be that the dual of a simplex should be a cell. Let's take this as the working definition of `manifold', and see where it takes us. Later we may see something  better, but this is a way to get started.\end{quotation}
We now know that the basic theory of PL manifolds is more elementary and accessible than either smooth or toplogical manifolds, and this working definition is a pretty good pointer to the theory. Why were Poincar\'e's successors slow to approach the subject this way, and when they did, why did it not work as well as we might have expected?

The first problem was that Poincar\'e and other nineteenth-century mathematicians objected to the use of explicit definitions. The objection goes  back  2400 years to Pythagoras and Plato, and is roughly that accepting a definition is like accepting a religious doctrine: you get locked in and blocked from any direct (intuitive) connection to `reality'. The precise-definition movement reflects practice in science: established definitions are distillations of the discoveries of our predecessors, and working definitions provide precise input needed for high-precision reasoning. It is odd that this aspect of scientific practice came so late to mathematics, but recall that in the nineteenth century there was still a strong linkage between mathematics and philosophy. And still to this day,  accepting a definition in philosophy is  like accepting a religious doctrine. 

An interesting transitional form appears  in a long essay by Tietze in 1908 (\cite{tietze}; see the translation at \cite{tietzeTrans}). He defined manifolds as polyhedra such that the link of a simplex is simply connected, but did not define `simply connected'. It is hard to imagine that he meant this literally. The use of the terminology `simply-connected' indicates familiarity with Poincar\'e's work with the fundamental group, but Poincar\'e  asks explicitly if it is possible for a 3-manifold ``to be simply-connected and yet not a sphere''. Simply-connected is obviously wrong one dimension higher.  His use of the term seems to have been a deliberately ambiguous placeholder in a  proposal for a ``big-picture'' view of manifolds. This reflects the philosophical idea that big pictures should be independent of details, and the goal of  heuristic arguments in the nineteenth-century tradition was to convince people that this was the right intuition, not actually prove things. On a practical level, 
Tieze may have been mindful of the advantages ambiguity had for Poincar\'e. People worked  hard trying to find interpretations of Poincar\'e's ideas that would make them correct, but could not have been so generous if he had tried to be more precise and guessed wrong. 

The next milestone we mention is the introduction of PL homology manifolds as a precise setting for the study of duality.  Wilder \cite{wilder} attributes this to Veblen in 1916. They knew these were not always locally euclidean so would not be the final context for geometric work, but  they would serve for algebraic topology until the geometric people got their acts together. 

In the geometric line at that time, people were experimenting with various precise replacements for Tieze's placeholder. The favorite was ``stars homeomorphic to Euclidean space''. Today we would see this as a mixed-category idea that for general reasons is unlikely to be correct and in any case is inappropriate for a basic definition. This experience was not available at the time, of course, but they were not having success with homeomorphisms and there were clues that an all-PL version would have advantages. Why did they stick with homeomorphisms for so long? There were two philosophical concerns and a technical problem.

The first philosophical concern was that a `manifold' should be a \emph{thing}. A topological space was considered a primitive thing\footnote{We now think of a topological space as a structure (a topology) on a set. In the Poincar\'e tradition, spaces were primitive objects with properties extrapolated from those of subsets of euclidean spaces.}, and a space that satisfies a property (eg.~locally homeomorphic to euclidean space) is also a thing. A simplicial complex is also a thing. A polyhedron, however, is a space with an equivalence class of triangulations. This is a \emph{structure} on a thing, not a primitive thing, so for philosophical reasons could not qualify as a correct definition of `manifold'. This objection also blocked the use of  coordinate charts to globalize differential structures.  

The second philosophical objection to PL manifolds has to do with the ``recognition problem''. A simplicial complex is a finite set of data. Suppose someone sent you one in the mail. How would you know whether or not links of simplices were PL equivalent to spheres? Suppose the sender enclosed a note asserting that this was so. How could you check to be sure it was true?  Bertrand Russell summarized the philosophical attitude toward such things \cite{russell}, p.~71
\begin{quote}
The method of ``postulating'' what we want has many advantages; they are the same as the advantages of theft over honest toil.
\end{quote}
The manly thing to do, then, is to \emph{prove}   links are PL spheres, and  \emph{assuming} this is cowardly and philosophically dishonest. Today we might wonder that assuming that a space is locally euclidean (rather than recognizing it as being so) is ok, while assuming PL is not.  At any rate one consequence was that  the generalized Poincar\'e conjecture\footnote{The generalized Poincar\'e conjecture is the assertion that a polyhedron that is known to be a PL manifold and that has the homotopy type of the sphere, is PL equivalent to the sphere.}  (then referred to as `the sphere problem') seemed to be essential to justify work in higher-dimensional PL manifolds.  The effect was to paralyze the field. 

The technical problem had to do with the definition of ``PL equivalence'' of simplicial complexes.   The modern definition is that they should have a common subdivision. This is very convenient technically because if you show some invariant does not change under a single subdivision then it must be a PL invariant. For traditionalists, however, it seemed uncomfortably existential. Equivalence of smooth or topological objects uses a nice concrete function with specific local properties; shouldn't PL follow this pattern? Brouwer, the great intuitionist, intuited a direct simplicial criterion for stars in simplicial complexes to be ``Euclidean'' and proposed this as a replacement for Tietze's placeholder. His intuiton was ineffective, however, and later shown to be wrong\footnote{In 1941, after the dust had settled, J.~H.~C.~Whitehead reviewed the various proposals from the 1920s. Brouwer's proposal was particularly dysfunctional, and one has to wonder if he had actually tried to work with it in any serious way.}.

We finally come to Kneser's triangulation questions. In 1924 he gave precise formulations of what would have to be done to show that  `polyhedron locally homeomorphic with euclidean space' really did give a theory as envisioned by Poincar\'e, Tietze, et al. Whether he intended it or not, one message was roughly ``enough sterile big-picture speculation; time to focus on what it would take to make it work.''   In particular, since
Poincar\'e's use of dual cells gives duality between homologies defined with two different triangulations, the uniqueness of triangulations  was needed to show Poincar\'e's claims about duality were correct.  It must have seemed scandalous that this was still unresolved a quarter-century after Poincar\'e made the claims. We might also see Hilbert's influence in the concise straight-to-the-point formulations. 

When Van der Vaerden surveyed manifold theory in 1928  he describe it as a ``battlefield of techniques''. There had been  advances in methodology but still no effective  definitions and big issues were still unsettled. In fact the situation was already improving.  In 1926 Newman \cite{newman26}  had published a version in which stars were still assumed homeomorphic to Euclidean space, but with complicated combinatorial conditions. This still didn't work, but in 1928 he published a revision  \cite{newman28} in which this was replaced by the common-subdivision version  still used in the mature theory.  PL topoogy was finally launched but, as it turned out, a bit too late. 

\subsection{Overtaken} Manifolds were supposed to be a setting for global questions in analysis, so smooth manifolds were the main goal. We have been following the PL topology developed to make sense of Poincar\'e's combinatorial ideas, but there are two other approaches that would have done this. The most effective is singular homology. This requires some algebraic machinery, but it is simpler than the PL development, much more general, and connects better with analytic use of sheaves, currents, and deRham cohomology. The second is less effective but closer to Poincar\'e's ideas: show that smooth manifolds have standard (piecewise-smooth) triangulations. Remember that Kneser called the uniqueness-of-triangulations question ``the Hauptvermutung'' (principal assertion) because it would show that simplicial homology is independent of the triangulation.  Either of the other approaches would have accomplished this, and therefore achieved the principal motivation of the PL development. The historical question should be: given the obvious importance of the questions, why did it take so long to find \emph{any} of these solutions\footnote{Existence of piecewise-smooth triangulations was shown in 1940 by Whitehead \cite{whitehead40}.}? Slow development in the PL track is only interesting because the others were slow as well. 

Singular homology probably developed slowly because it is so far outside the received wisdom from Poincar\'e. It requires algebraic apparatus and while we now see plenty of clues about this in Poincar\'e's work, these only became visible after Noether's promotion of abstract algebra as a context for such things. \v Cech's open-cover approach to homology also pushed things in this direction, but again this was outside Poincar\'e's vision. 

The lack of an effective definition held up development of smooth manifolds, just as it held up development of PL. And, like PL and unlike the singular theory, no huge technical leaps were necessary: the main obstructions were ineffective intuitions and philosophical objections to structures defined with coordinate charts. These were finally overcome by  
 Veben and Whitehead \cite{vebWhitehead32} in 1932, and we can identify two things that made the advance possible. The first was a change to a more modern style that better reflects mathematical structure. Veblen and Whitehead did not give a philosophical argument or a speculative `big picture'; they developed enough basic structure (with technical details) to demonstrate conclusively that this was an effective setting for differential geometry. The second change was in the mathematical community. Young people were attracted by the power and depth of precise definitions and full-precision reasoning, and were more than ready to trade philosophy for success, while the old people committed to philosophy were fading away. These changes led to a great flowering of the differential theory, and it was the setting for some of the deepest and most remarkable discoveries of the second half of the twentieth century. 
 
One consequence of the smooth-manifold flowering, and the development of singular homology, was  a near abandonment of  PL topology for several decades. It continued to be used in low dimensions due to low-dimensional simplifications (homology identifies 2-manifolds). Enough of a community had been established to sustain some general activity, but it lacked the guidance of an important goal. 

The 1950s and 60s saw a renascence in PL topology. Smale's development of handlebody theory, and particularly his proof of a form of the generalized Poincar\'e conjecture, electrified the manifold communities. Smale's work was in the smooth world, coming from a study of the dynamics of Morse functions, but handles appear much more easily and naturally in PL.  Milnor's discovery of multiple smooth structues on the 7-sphere \cite{milnor56} was a huge boost. The reason was that Smale had proved that high-dimensional \emph{smooth} homopy spheres were \emph{homeomorphic} to the sphere. Stallings then used PL techniques to show that  high-dimensional PL homopy spheres were homeomorphic to the sphere. Both of these had the defect that the conclusions were in a different category from the hypotheses. Smale improved on this by using PL versions of his techniques to show that high-dimensional PL homtopy spheres are PL isomorphic to spheres. Milnor's discovery showed that this is false in the smooth world, so PL is genuinely simpler and closer to the original intuitions\footnote{Milnor's discovery also invalidated the intuition, inherited from Poincar\'e, that there would be a single world of `manifolds' where all techniques would be available. Subsequent developments, as we have seen here, revealed how confining that intuition had been.}. All this took place in the modern links-are-PL-spheres context. Kneser's questions played no role and, as far as the main-line developments were concerned, were a dead-end curiosity. 

 By the end of the 60s PL was again overshadowed, this time by development of purely topological manifold theory. Basic topological techniques are much more complicated than PL, almost insanely so in some cases, but the outcomes are more systematic and coherent.   Further progress on what seemed to be PL issues  also required outside techniques:  The 3-dimensional Poincar\'e conjecture was settled by Perelman with delicate analytic arguments almost 80 years after Kneser's work, and 100 years after Poincar\'e hinted that this might be the key to further progress.  The 4-dimensional case is still open in 2013, and no resolution is in sight.  Finally, as we have seen here, insight into the structure of homology 3-spheres seems to require gauge theory. 

\subsection{Summary}
 Kneser's triangulation questions were a careful formulation of what it would take to develop a theory of manifolds that followed Poincar\'e's intuitions and nineteenth-century philosophy. Not long after, more fruitful approaches emerged based on  full-precision  twentieth-century methodology. Kneser's questions proved to be a curiosity: an nice challenge for developing technology, but apparently without significant implications. 

\bibliographystyle{plain}
\cleardoublepage

\bibliography{triangexpo}

\end{document}